\documentclass{amsart}%{imanum}
\usepackage{natbib}

\numberwithin{equation}{section}
\usepackage{hyperref}
\swapnumbers
\theoremstyle{plain}
\newtheorem{theorem}[subsection]{Theorem}%[section]
\newtheorem{lemma}[subsection]{Lemma}
\newtheorem{corollary}[subsection]{Corollary}
\newtheorem{proposition}[subsection]{Proposition}
\theoremstyle{definition}
\newtheorem{definition}[subsection]{Definition}

\newtheorem{remark}[subsection]{Remark}
\providecommand{\linkedemail}[1]{\href{mailto:#1}{#1}}

\newcommand{\be}{\begin{equation}}
\newcommand{\ee}{\end{equation}}
\newcommand{\bea}{\begin{eqnarray}}
\newcommand{\eea}{\end{eqnarray}}
\newcommand{\beaa}{\begin{eqnarray*}}
\newcommand{\eeaa}{\end{eqnarray*}}
\newcommand{\suint}{\sum_{j=1}^{m-1}\int_{t^{j-1}}^{t^j}}
\newcommand{\mtau}{\int_{t^{m-1}}^{\tau}}

\newcommand{\norm}[2]{\ensuremath{\| #1 \|_{#2}}}

\newcommand{\mbf}[1]{\mbox{\boldmath$\rm{#1}$}}

\newcommand{\el}{ \kappa \in \mathcal{T}}

\DeclareMathOperator{\diam}{diam}

\newcommand{\mean}[1]{\{ \kern -1.6mm \{#1\} \kern -1.6mm \}}
\newcommand{\ha}{\frac{1}{2}}
\newcommand{\jump}[1]{[ \kern -.7mm [#1] \kern -.7mm ]}

\newcommand{\ud}{\mathrm{d}}
\newcommand{\ndg}[1]{| \kern -.25mm \|{#1}| \kern -.25mm \|}
\newcommand{\ltwo}[2]{\|{#1}\|_{#2}}
\newcommand{\enorm}[1]{\|{#1}\|_a}
\newcommand{\linf}[2]{\|{#1}\|_{L^{\infty}({#2})}}

\newcommand{\mesh}[1]{\mathcal{T}^{{#1}}}

\newcommand{\Uh}{U}%%{\hat{U}}
\newcommand{\wh}{w}%%{\hat{w}}

\DeclareMathOperator{\esssup}{ess\sup}
\newcommand{\Const}[1]{C_{\mathrm{#1}}}

\begin{document}

\title[$L^{\infty}(L^2)$-norm a posteriori bounds for wave problems]{A posteriori $L^{\infty}(L^2)$-error bounds for finite element approximations to the wave equation}
%%\shorttitle{$L^{\infty}(L^2)$-norm a posteriori bounds for wave problems}

%%\author
{%
\author{Emmanuil H.~Georgoulis}
\address{
Emmanuil H.~Georgoulis
\\
Department of Mathematics
\\
University of Leicester
\\
University Road
\\
Leicester
\\
LE1~7RH
\\
United Kingdom.
}
\email{\linkedemail{Emmanuil.Georgoulis@le.ac.uk}}
\author{Omar Lakkis}
\address{
  Omar Lakkis
  \\
  Department of Mathematics
  \\
  University of Sussex
  \\
  Brighton
  \\
  GB-BN1~9RF
  \\
  England UK.
}
\email{\linkedemail{O.Lakkis@sussex.ac.uk}}
\author{Charalambos Makridakis}
\address{
Department of Applied Mathematics,
\\
University of Crete, GR-71409 Heraklion, Greece, 
\\
and
\\
Institute for Applied and Computational Mathematics,
\\
Foundation for Research and Technology-Hellas,
\\
Vasilika Vouton P.O. Box 1527, Heraklion
\\
GR-71110
\\
Greece.
}
\email{\linkedemail{makr@tem.uoc.gr}}
% Short list of authors for running heads:
%%\shortauthorlist{E.H.~Georgoulis, O.~Lakkis, C.~Makridakis}

\maketitle

\begin{abstract}
{We address the error control of
Galerkin discretisation (in space) of linear second order hyperbolic problems.
More specifically, we derive a posteriori error bounds in the $L^{\infty}(L^2)$-norm for finite element methods
for the linear wave equation, under minimal regularity assumptions.
The theory is developed for both the space-discrete case, as well as
for an implicit fully discrete scheme. The derivation of these bounds relies crucially on
carefully constructed space- and time-reconstructions of the discrete numerical solutions, in conjunction with
a technique introduced by Baker  (SIAM J. Numer. Anal., 13, 1976)
in the context of a priori error analysis of Galerkin discretisation
of the wave problem in weaker-than-energy spatial norms.}
{a posteriori bounds, wave equation, implicit time-stepping, reconstruction}
\end{abstract}

%\end{document}

\section{Introduction}
In computing approximate solutions of evolution initial-boundary
value problems mesh-adaptivity plays an important role, in that it
drives variable resolution requirements, thereby contributing
reduction in computational cost. Adaptive strategies are often based
on a posteriori error estimates, i.e., computable quantities which
estimate the error of the finite element method measured in a
suitable norm (or other functionals of interest).

%\end{document}

A posteriori error bounds are well developed for stationary boundary
value problems \citep[e.g.,][and
the references therein]{verfurth, ainsworth_oden,
  babuska_strouboulis, cc-zz, dorfler, stevenson, cascon_nochetto}. Adaptivity and error estimation for
parabolic problems has also been an active area of research for the
last two decades \citep[e.g.,][and
the references therein]{eriksson-johnson:2,
  verfurth_parabolic,picasso,houston-suli_lagrange,makr-noch,bergam-bernardi-mghazli:05,
  bernardi-verfuerth:05,lakkis-makr}.

Surprisingly, there has been considerably less work on the error
control of finite element methods for second order hyperbolic
problems, despite the substantial amount of research in the design
of finite element methods for the wave problem  
\citep[e.g.,][and the references therein]{baker_wave,baker_bramble,baker-dougalis-serbin,baker-dougalis2, dougalis-serbin2,johnson_wave,makr1,joly,joly2,joly-tsogka2,kar_makr}. A posteriori bounds for standard
implicit time-stepping finite element approximations to the linear
wave equation have been proposed and analyzed (but only in very
specific situations) by \citet{adjerid}.  Also, \citet{bernardi_suli} derive rigorous a posteriori
bounds, using energy arguments, for finite element methods with first order implicit time-stepping. Moreover, \citet{bernardi_suli} proposes an adaptive algorithm based on the a posteriori bounds derived therein.  Goal-oriented error estimation for wave problems (via duality
techniques) is also available~\citet{bangerth-rannacher:99,bangerth-rannacher:01}, while some
earlier work on a posteriori estimates for first order hyperbolic
systems have been studied in the time semidiscrete setting~\citet{makr_noch_diss}, as well as in the fully discrete one~\citet{johnson_wave,suli,suli2}.

In this work, we derive a posteriori bounds in the
$L^{\infty}(L^2)$-norm of the error.  The theory is developed for both the
space-discrete case, as well as for the practically relevant case of
an implicit fully discrete scheme. The derivation of these bounds relies crucially on
\emph{reconstruction} techniques, used earlier for parabolic
problems~\citet{makr-noch,lakkis-makr,amn}. Another key tool in our analysis is the special testing procedure due to \citet{baker_wave}, who used it in the a priori error analysis
of Galerkin discretisation of the wave problem in weaker-than-energy
spatial norms. It is expected (although it is not considered here) that the novel space-time-reconstruction presented here could also be applicable to spatially non-linear second order hyperbolic problems and to different, possibly non-confirming, spatial discretisations. Moreover, it is also possible to combine the abstract results presented here with a wide class of a posteriori error estimators for elliptic problems. 

While for the proof of a posteriori bounds for the semidiscrete
case, the \emph{elliptic reconstruction} previously considered
in~\citet{makr-noch,lakkis-makr} suffices, the fully discrete
analysis necessitates the careful introduction of a novel space-time
reconstruction, satisfying a crucial \emph{local vanishing moment
property} in time.  Our approach is based on the one-field
formulation of the wave equation and, thus, non-trivial three-point
time reconstructions are required. A further challenge presented by
the wave equation is the special treatment of deriving bounds for
the ``elliptic error'' of the reconstruction framework, to obtain
practically implementable residual estimators. The derived a
posteriori estimators are formally of optimal order, i.e., of the
same order as the error on uniform space- and time-meshes.

The a posteriori bounds proposed in this work could be used within an adaptive algorithm, such as the one presented in \citet{bernardi_suli}. However, this is a important task in its own right and will be considered elsewhere.

The rest of this work is organized as follows. In \S\ref{prelim} we
present the model problem and the necessary basic definitions along
with the finite element methods for the wave equations considered in
this work. In \S\ref{abs_apost} we consider the case of a posteriori
bounds for the space-discrete problem.  In \S\ref{full_discrete}, we
derive abstract a posteriori error bounds for the fully-discrete
implicit finite element method, while in
\S\ref{sec:residual-estimates} the case of a posteriori bounds of
residual type are presented. In  \S\ref{conclusions}, we draw some
final concluding remarks.

\section{Preliminaries}\label{prelim}
\subsection{Model problem and notation}\label{mod}
We denote by $L^p(\omega)$, $1\le p\le +\infty$, for
$\omega\subset\mathbb{R}^d$ open, the Lebesgue spaces, with corresponding
norms $\|\cdot\|_{L^p(\omega)}$. The norm of $L^2(\omega)$, denoted
by $\ltwo{\cdot}{\omega}$, corresponds to the $L^2(\omega)$-inner
product $\langle \cdot,\cdot\rangle_{\omega}$. We denote by
$H^s(\omega)$, the Hilbertian Sobolev space of order $s\ge0$ of
real-valued functions defined on $\omega\subset\mathbb{R}^d$ (see, e.g., \citet{adams} for definitions and basic properties); in
particular $H^1_0(\omega)$ signifies the space of functions in
$H^1(\omega)$ that vanish on the boundary $\partial\omega$ (boundary
values are taken in the sense of traces). Negative order Sobolev
spaces $H^{-s}(\omega)$, for $s>0$, are defined through duality. In
the case $s=1$, the definition of $\langle
\cdot,\cdot\rangle_{\omega}$ is extended to the standard duality
pairing between $H^{-1}(\omega)$ and $H^1_0(\omega)$.  For $1\le
p\le +\infty$, we also define the spaces $L^p(0,T,X)$, with $X$
being a real separable Banach space with norm $\|\cdot\|_X$, consisting of all
measurable functions $v: (0,T)\to X$, for which
\begin{equation}
  \begin{aligned}
    \|v\|_{L^p(0,T;X)}&:=\Big(\int_0^T \|v(t)\|_{X}^p\ud t\Big)^{1/p}<+\infty,\quad \text{for}\quad 1\le p< +\infty, \\
    \|v\|_{L^{\infty}(0,T;X)}&:=\esssup_{0\le t\le T}\|v(t)\|_{X}<+\infty,\quad \text{for}\quad  p = +\infty.
  \end{aligned}
\end{equation}
%Finally, we denote by $C^m(0,T;X)$ the space of functions $v: (0,T)\to X$ with $m$ continuous derivatives
%with respect to time, with norm
%\begin{equation}\|v\|_{C^m(0,T;X)}:=\max_{0\le k\le m}\sup_{0< t<T}\|\frac{\partial^k}{\partial t^k}v(t)\|_{X}<+\infty.\end{equation}

Let $\Omega\subset \mathbb{R}^d$ be a bounded open polygonal domain with
Lipschitz boundary $\partial\Omega$. For brevity, the standard inner product on $L^2(\Omega)$ will be
denoted by $\langle \cdot,\cdot\rangle$ and the corresponding norm by $\ltwo{\cdot}{}$.

For time $t\in(0,T]$, we consider the linear second order hyperbolic initial-boundary value problem of
finding $u\in L^2(0,T;H^1_0(\Omega))$, with $u_t\in L^2(0,T;L^2(\Omega))$ and $u_{tt}\in L^2(0,T;H^{-1}(\Omega))$
such that
\begin{equation}\label{pde}
 u_{tt}-\nabla \cdot (a\nabla u)=f\quad\text{in } (0,T)\times\Omega,
\end{equation}
where $f\in L^2(0,T; L^2(\Omega))$ and
$a$ is a scalar-value function in $\in C(\bar{\Omega})$, with $0<\alpha_{\rm min}\le a\le \alpha_{\rm max}$, such that
\begin{equation}
  \label{modelbc}
  \begin{split}
    u(x,0)&=u_0(x)\text{ on }\Omega\times\{0\},\\
    u_t(x,0)&=u_1(x)\text{ on }\Omega\times\{0\}\\
    u(0,t)&= 0\text{ on }\partial\Omega\times(0,T],
  \end{split}
\end{equation}
where $u_0\in H^1_0(\Omega)$ and $u_1\in L^2(\Omega)$. For existence and uniqueness results for this problem, we refer, e.g., to \citet{raviart}, Chapter 8.

%%In standard fashion,
We identify
a function $v\in\Omega\times [0,T]\to \mathbb{R}$ with the function $v:[0,T]\to H^1_0(\Omega)$ and
we use the shorthand $v(t)$ to indicate $v(\cdot,t)$.

\subsection{Finite element method}
Let $\mathcal{T}$ be a shape-regular subdivision of $\Omega$ into
disjoint open simplicial or quadrilateral elements. Each element $\el$ is constructed via mappings
$F_{\kappa}:\hat{\kappa}\to\kappa$, where
$\hat{\kappa}$ is the reference simplex or reference square, so that
$\bar{\Omega}=\cup_{\el}\bar{\kappa}$, \citep[see, e.g.,][]{ciarlet}.

For a nonnegative integer $p$, we denote by
$\mathcal{P}_p(\hat{\kappa})$ either the set of all polynomials on
$\hat{\kappa}$ of degree $p$ or less, when $\hat{\kappa}$ is the
simplex, or the set of polynomials of at most degree $p$ in each
variable, when $\hat{\kappa}$ is the reference square (or cube). We
consider $p$ fixed and use the finite element space
\begin{equation}\label{eq:FEM-spc}
V_h :=\{v\in H^1_0(\Omega):v|_{\kappa}\circ F_{\kappa}
    \in \mathcal{P}_p(\hat{\kappa}),\,\el\}.
\end{equation}

Further, we denote by $\Gamma:=\cup_{\kappa\in\mathcal{T}}(\partial\kappa\backslash\partial\Omega)$, i.e.,
the union of all $(d-1)$-dimensional element edges (or faces) $e$ in $\Omega$
associated with the subdivision $\mathcal{T}$ excluding the
boundary. We introduce the mesh-size function $h:\Omega\to \mathbb{R}$, defined
by $h(x)= \diam{\kappa}$, if $x\in
\kappa$ and $h(x)= \diam(e)$, if $x\in e$ when $e$ is an edge.

The semidiscrete finite element method for the
initial-boundary value problem (\ref{pde})--(\ref{modelbc}) consists in
finding $U\in L^2(0,T;V_h)$ such that
\begin{equation}\label{fem_semi}
\langle U_{tt},V\rangle+a(U,V)=\langle f,V\rangle
\quad \forall V\in L^2(0,T;V_h),
\end{equation}
where the bilinear form $a$ is defined for each $z,v\in H^1_0(\Omega)$ by
\begin{equation}
a(z,v)=\int_{\Omega} a\nabla z\cdot\nabla v\,\ud x,
\end{equation}
and the corresponding energy norm is defined for $v\in H^1_0(\Omega)$ by
\begin{equation}\label{enorm}
\enorm{v}=\ltwo{\sqrt{a}\nabla v}{}.
\end{equation}

To introduce the fully-discrete implicit scheme approximating
(\ref{pde})--(\ref{modelbc}), we consider a subdivision of the time interval $(0,T]$ into
subintervals $(t^{n-1},t^{n}]$, $n=1,\dots, N$,
with $t^0=0$ and $t^N=T$, and we define $k_n:=t^{n}-t^{n-1}$, the local time-step.
Associated with the time-subdivision, let $\mesh{n}$, $n=0,\dots, N$, be a sequence of meshes which are
assumed to be \emph{compatible} (see, e.g., \cite{lakkis-makr} for a precise definition of mesh compatibility in this context), in the sense that for any two
consecutive meshes $\mesh{n-1}$ and $\mesh{n}$, $\mesh{n}$ can be obtained
from $\mesh{n-1}$ by locally coarsening some of its elements and then locally refining some (possibly other) elements.
The finite element space corresponding to $\mesh{n}$ will be denoted by $V_h^n$.

We consider the fully discrete scheme for the wave problem (\ref{pde}), (\ref{modelbc})
\begin{equation}\label{fem_fd}
  \begin{split}
    &\text{for each $n=1,\dots,N$, find}\ U^{n}\in V_h^{n}\text{ such that  }
    \\
    &\langle \partial^2 U^{n},V\rangle+a(U^{n},V)=\langle f^{n},V\rangle
    \quad \forall V\in V_h^{n},
  \end{split}
\end{equation}
where $f^{n}:=f(t^{n},\cdot)$, the backward second and first finite differences
\begin{equation}\label{second_der_def}
\partial^2 U^{n}:=\frac{\partial U^{n}-\partial U^{n-1}}{k_n},\quad
\end{equation}
with
\begin{equation}\label{second_der_first_def}
\partial U^{n}:=
\begin{cases}
  \displaystyle\frac{U^{n}- U^{n-1}}{k_{n}},
  &\text{ for }n=1,2,\dots,N,
  \\
  V^0:=\pi^0 u_1
  &\text{ for }n=0,
\end{cases}
\end{equation}
%% i.e., , respectively. The first step, for $n=1$,
%% is defined by
%% \begin{equation}\label{fem_fd_neq1}
%% \text{find}\ U^{1}\in V_h^{1}\text{ such that  }
%% \langle \partial^2 U^{1},V\rangle+a(U^{1},V)=\langle f^{1},V\rangle
%% \quad \forall V\in V_h^{1},
%% \end{equation}
%% where $f^{1}:=f(t^{1},\cdot)$,
where $U^0:=\pi^0 u_0$, and $\pi^0: L^2(\Omega)\to V_h^0$ a suitable projection
onto the finite element space (e.g., the orthogonal $L^2$-projection operator). 

Denoting by $A^n$ the stiffness matrix for the mesh $\mesh{n}$, and by $\underline{U}^n$ the respective coefficient vector for $U^n$, the implicit method reads: find $\underline{U}^n\in \mathbb{R}^{\dim V^n_h}$ such that
\[
\frac{1}{k_n}\Big(\frac{\underline{U}^n-\underline{U}^{n-1}}{k_n}-\frac{\underline{U}^{n-1}-\underline{U}^{n-2}}{k_{n-1}}\Big) + A^n \underline{U}^n = \underline f^{n},
\]
with $\underline{f}^n:=(f_i^n)_{i}$, with $f_i^n:=(f^n,\phi_i^n)$, $i=1,\dots, \dim V^n_h$, for $\phi_i^n$ such that $V_h^n={\rm span}\{\phi_i^n:i=1,\dots, \dim V^n_h \}$.

\section{A posteriori error bounds for the semi-discrete problem}\label{abs_apost}

We derive here a posteriori error bound for the error
$\linf{u-U}{0,T;L^2(\Omega)}$ between the exact solution
of~(\ref{pde}), (\ref{modelbc}) and that of the semidiscrete scheme
\ref{fem_semi}.
%We first derive an abstract result in
%Theorem~\ref{apost_semi_thm_wave_l2}, which independent of which
%elliptic estimators one plans to use.  We then recall residual a
%posteriori estimates for elliptic equation
%in~\ref{mraposterioritheorem} and apply it to derive a ``concrete''
%result for the wave equation in
%Corollary~\ref{apost_semi_thm_wave_l2_cor}.

\begin{definition}[Elliptic reconstruction and error splitting]\label{elliptic_recon}
Let $U$ be the (semidiscrete) finite element solution to the problem (\ref{fem_semi}).
Let also $\Pi:L^2(\Omega)\to V_h$ be the orthogonal $L^2$-projection operator onto the finite element space
$V_h$.
We define the \emph{elliptic reconstruction} $w=w(t)\in H^1_0(\Omega)$, $t\in [0,T]$, of $U$ to be
the solution of the elliptic problem
\begin{equation}\label{ell_rec}
a(w, v) = \langle g, v\rangle\quad \forall v\in H^1_0(\Omega)%%,\ \text{a.e.}\ t\in [0,T],
\end{equation}
where
\begin{equation}\label{tzi}
g:= A U-\Pi f +f,
\end{equation}
and $A:V_h\to V_h$ is the \emph{discrete elliptic operator} defined by
\begin{equation}\label{ell_rec2}
 \langle A q, \chi \rangle = a(q,\chi)\quad \forall q,\chi\in V_h.
\end{equation}
We decompose the error as follows:
\begin{equation}\label{splitting_semi_wave}
e:=U-u=\rho-\epsilon,\ \text{where}\ \epsilon:=w-U,\ \text{and}\ \rho:=w-u.
\end{equation}
\end{definition}

\begin{lemma}[Error relation]\label{energy_arg_wave}
With reference to the notation in (\ref{splitting_semi_wave}) we have
\begin{equation}\label{en_ar_wave}
\langle e_{tt},v\rangle +a(\rho,v)=0\quad\forall v\in H^1_0(\Omega).
\end{equation}
\end{lemma}
\begin{proof}  We have, respectively,
\begin{equation}\label{rho-term_wave}
\begin{aligned}
\langle  e_{tt},v\rangle +a(\rho,v)
&=\langle  U_{tt}, v   \rangle +a(w,v)-\langle  u_{tt}, v   \rangle -a(u,v) \\
&=\langle  U_{tt}, v   \rangle +a(w,v) - \langle f,v\rangle\\
&=\langle  U_{tt},\Pi v\rangle +a(w,v) - \langle f,v\rangle\\
&= - a(U,\Pi v) +a(w,v)+\langle  \Pi f-f,v\rangle =0,
\end{aligned}
\end{equation}
observing the identity $a(U,\Pi v)-\langle  \Pi f-f,v\rangle= a(w,v)$ due to the construction of $w$.\end{proof}

\begin{theorem}[Abstract semidiscrete error bound]\label{apost_semi_thm_wave_l2}
With the notation introduced in (\ref{splitting_semi_wave}), the
following error bound holds:
\begin{equation}\label{apost_semi_wave}
\begin{aligned}
\linf{e}{0,T;L^2(\Omega)}\le& \linf{\epsilon}{0,T;L^2(\Omega)}+\sqrt{2}\Big(\ltwo{u_0-U(0)}{}
+\ltwo{\epsilon(0)}{}\Big)\\
&+2\int_0^T\ltwo{\epsilon_t}{} + C_{a,T} \ltwo{u_1- U_t(0)}{},
\end{aligned}
\end{equation}
where $C_{a,T}:=\min\{2T,\sqrt{2C_{\Omega}/\alpha_{\rm min}}\}$, where $C_{\Omega}$ is
the constant of the Poincar\'e--Friedrichs inequality
$\ltwo{v}{}^2\le C_{\Omega}\ltwo{\nabla v}{}^2$, for $v\in H^1_0(\Omega)$.
\end{theorem}
\begin{proof}
We use a testing procedure due to \citet{baker_wave}.
Let $\tilde{v}:[0,T]\times\Omega\to\mathbb{R}$ with
\begin{equation}\label{tilde_v}
\tilde{v}(t,\cdot)=\int_t^{\tau} \rho(s,\cdot)\ud s,\quad t\in[0,T],
\end{equation}
from some fixed $\tau\in[0,T]$. Clearly $\tilde{v}\in H^1_0(\Omega)$ as $\rho\in H^1_0(\Omega)$. Also, we observe that:
\begin{equation}\label{properties_tilde_v}
\tilde{v}(\tau,\cdot)=0,\quad \nabla\tilde{v}(\tau,\cdot)=0,\quad\text{and}\quad
\tilde{v}_t(t,\cdot)=-\rho(t,\cdot),\quad \text{a.e. in}\ [0,T].
\end{equation}
Set $v=\tilde{v}$ in (\ref{en_ar_wave}), integrate between $0$ and $\tau$ with respect
to the variable $t$ and integrate by parts the first term on the left-hand side, to obtain
\begin{equation}
- \int_0^{\tau} \langle  e_t,\tilde{v}_t\rangle
+ \langle e_t(\tau),\tilde{v}(\tau)\rangle
- \langle  e_t( 0  ),\tilde{v}( 0  )\rangle
+ \int_0^{\tau} a(\rho,\tilde{v})
= 0.
\end{equation}
Using (\ref{properties_tilde_v}), we have
\begin{equation}
 \int_0^{\tau}\ha \frac{d}{dt}\ltwo{\rho(t)}{}^2
- \int_0^{\tau} \ha \frac{d}{dt}a(\tilde{v}(t),\tilde{v}(t))
= \int_0^{\tau}\langle \epsilon_t,\rho\rangle
+ \langle e_t( 0  ),\tilde{v}( 0  )\rangle,
\end{equation}
which implies
\begin{equation}
 \ha\ltwo{\rho(\tau)}{}^2 - \ha\ltwo{\rho(0)}{}^2
+\ha a(\tilde{v}(0),\tilde{v}(0))
= \int_0^{\tau}\langle \epsilon_t,\rho\rangle
+ \langle e_t( 0  ),\tilde{v}( 0  )\rangle.
\end{equation}
Hence, we deduce
\begin{equation}\label{abstract_bound_semi_mid}
 \ha\ltwo{\rho(\tau)}{}^2 - \ha\ltwo{\rho(0)}{}^2+\ha a(\tilde{v}(0),\tilde{v}(0))
\le \max_{0\le t\le T} \ltwo{\rho(t)}{} \int_0^{\tau}\ltwo{\epsilon_t}{}
+ \ltwo{e_t( 0  )}{}\ltwo{\tilde{v}( 0  )}{}.
\end{equation}
Now, we select $\tau$ such that
$\ltwo{\rho(\tau)}{}=\max_{0\le t\le T}\ltwo{\rho(t)}{}$, (this is possible due to the continuity of $u$ in the time variable under the data and domain regularity assumptions above, see, e.g., \citet{raviart}, Chapter 8,)  and we
present two alternative, but complementary, ways to complete the
proof.

In the first way, we start by observing that, for this $\tau$, we have
$\ltwo{\tilde{v}(0)}{}\le \ltwo{\rho(\tau)}{}$, which gives
\begin{equation}\label{intermediate_semi_dis}
 \frac{1}{4}\ltwo{\rho(\tau)}{}^2 - \ha\ltwo{\rho(0)}{}^2
\le \Big(\int_0^{\tau}\ltwo{\partial_t \epsilon}{}
+ \tau\ltwo{e_t( 0  )}{}\Big)^2.
\end{equation}
Using the bound $\ltwo{\rho(0)}{}\le \ltwo{e(0)}{}+\ltwo{\epsilon(0)}{}$,
$e(0)=U(0)-u_0$ and $e_t(0)=U_t(0)-u_1$, and (\ref{intermediate_semi_dis}) for $\tau$ as above, we conclude that
\begin{equation}\label{apost_semi_wave_first}
\begin{aligned}
\linf{e}{0,T;L^2(\Omega)}\le &\linf{\epsilon}{0,T;L^2(\Omega)}+\linf{\rho}{0,T;L^2(\Omega)} \\
\le & \linf{\epsilon}{0,T;L^2(\Omega)}+\sqrt{2}\Big(\ltwo{u_0-U(0)}{}
+\ltwo{\epsilon(0)}{}\Big)\\
&+2\Big(\int_0^T\ltwo{\epsilon_t}{} + T \ltwo{u_1-U_t(0)}{}\Big).
\end{aligned}
\end{equation}

The second alternative, described next, consists in a different
treatment of the last term on the right-hand side of
(\ref{abstract_bound_semi_mid}). The Poincar\'e--Friedrichs inequality
and the positivity of the diffusion coefficient $a$ imply
$\ltwo{\tilde{v}( 0 )}{}^2\le C_{\Omega}\alpha_{\rm
  min}^{-1}\ltwo{\tilde{v}( 0 )}{a}^2$, for some constant $C_{\Omega}$
depending on the domain $\Omega$ only. Combining this bound with
(\ref{abstract_bound_semi_mid}), we arrive to
\begin{equation}\label{apost_semi_wave_almost_second}
 \ha\ltwo{\rho(\tau)}{}^2 - \ha\ltwo{\rho(0)}{}^2
\le \max_{0\le t\le T} \ltwo{\rho(t)}{} \int_0^{\tau}\ltwo{\epsilon_t}{}
+ \ha C_{\Omega}\alpha_{\rm min}^{-1}\ltwo{e_t( 0  )}{}^2,
\end{equation}
which implies
\begin{equation}\label{apost_semi_wave_second}
\begin{aligned}
\linf{e}{0,T;L^2(\Omega)}\le& \linf{\epsilon}{0,T;L^2(\Omega)}+\sqrt{2}\Big(\ltwo{u_0-U(0)}{}
+\ltwo{\epsilon(0)}{}\Big)\\
&+2\int_0^T\ltwo{\epsilon_t}{} + \sqrt{2C_{\Omega}/\alpha_{\rm min}} \ltwo{u_1-U_t(0)}{}.
\end{aligned}
\end{equation}
Taking the minimum of the bounds (\ref{apost_semi_wave_first}) and (\ref{apost_semi_wave_second}) yields the
result.
\end{proof}

\begin{remark}[Short and long integration times]
The use of two alternative arguments in the last step of the proof of
Lemma \ref{energy_arg_wave} improves the ``reliability
constant'' $C_{a,T}$ that works for both the short-time and the
long-time integration regimes.
\end{remark}

\begin{remark}[Completing the a posteriori estimation]
  To obtain a practical a posteriori bound, we need to estimate the norms involving the elliptic error $\epsilon$.
  By construction, the elliptic reconstruction $w$
  is the exact solution to the elliptic boundary-value problem (\ref{ell_rec}) whose
  finite element solution is $U$. Indeed,
  inserting $v=V\in V_h$ in (\ref{ell_rec}), we have
  \begin{equation}\label{crucial_ell_rec_property}
    a(w,V)=\langle AU-\Pi f +f,V\rangle=a(U,V),
  \end{equation}
  which
  implies the Galerkin orthogonality property $a(w-U,V)=0$. Therefore,
  by construction, $\epsilon$ is the error of the finite element method on $V_h$ for the
  elliptic problem
  \begin{equation}
  -\nabla\cdot(a\nabla w)= g,
  \end{equation}
  with homogeneous Dirichlet boundary conditions, with $g$ defined by (\ref{tzi}).
\end{remark}

\begin{definition}
  For every element face $e\subset\Gamma$, we define the \emph{jump} across $e$ of
  a field $\mbf w$, defined in an open neighborhood of $e$, by
  \begin{equation}
    \jump{\mbf w}(x)=\lim_{\delta\to 0}\big(\mbf w(x+\delta\mbf{n}_e)-\mbf w(x-\delta\mbf{n}_e)\Big)\cdot\mbf n_e,
  \end{equation}
  for $x\in e$, where $\mbf{n}_e$ denotes one of the two normal vectors
  to $e$ (the definition of jump is independent of the choice).
\end{definition}

\begin{theorem}[Elliptic a posteriori residual bounds]
\label{mraposterioritheorem}
Let $z \in H^1_0(\Omega)$ be the solution to the elliptic problem:
\begin{equation}\label{ell_prob}
-\nabla\cdot(a\nabla z)=r
\end{equation}
$r\in L^2(\Omega)$ and $\Omega$ convex, and let $Z\in V_h$ be the
finite element approximation of $z$ satisfying
\begin{equation}
  a(Z,V)=\langle r,V\rangle\quad\forall V\in V_h.
\end{equation}
Then, there exists a positive constant $\Const{el}$, independent of
$\mathcal{T}$, $h$, $z$ and $Z$, so that
\begin{equation}\label{apost_bound_bound}
\norm{z-Z}{}^2 \le \Const{el}\, \mathcal{E}(Z,r,\mathcal{T}),
\end{equation}
where
\begin{equation}
  \mathcal{E}(Z,r,\mathcal{T}) := \Big(\sum_{\kappa \in
    \mathcal{T}} \Big( \norm{h^2(r+\nabla\cdot(a\nabla
    Z)}{\kappa}^2 + \sum_{e \subset\Gamma}
  \norm{h^{3/2}\jump{a\nabla Z}}{e}^2\Big)\Big)^{1/2}.
\end{equation}
\end{theorem}
Such results (some with various extra assumptions) are generally available in the literature. We refer to  {\citet[Remark 2.4]{verfurth}, \citet[Theorem 2.7]{ainsworth_oden}} for proofs of Theorem \ref{ell_prob}, and to the references therein for similar approaches. 

%\textcolor{red}{ The following result is very easy. It will be better to just mention it--not as theorem--
%when it will be used}For further reference, we also note, without proof,
%a (standard) stability result for the elliptic boundary-value problem.
%
%\begin{theorem} \label{pde_stability_theorem}
%Let $z,z_* \in H^1_0(\Omega)$ be the solutions to the elliptic problems:
%\begin{equation}\label{ell_prob2}
%-\nabla\cdot(a\nabla z)= r\quad\text{and}\quad -\nabla\cdot(a\nabla z_*)= r_*
%\end{equation}
%with homogeneous Dirichlet boundary conditions, respectively. Then, we have
%\begin{equation}\label{pde_stability}
%\norm{z-z_*}{} \le \Const{\Omega}\alpha_{\min}^{-1}\norm{r-r_*}{}.
%\end{equation}
%\end{theorem}
%

\begin{corollary}[Semidiscrete residual-type a posteriori error bound]\label{apost_semi_thm_wave_l2_cor}
Assume that the hypotheses of Theorems \ref{apost_semi_thm_wave_l2} and \ref{mraposterioritheorem} hold.
Assume further that $f$ is differentiable with respect to time. Then the following error bound holds:
\begin{equation}\label{apost_semi_wave_cor}
\begin{aligned}
\linf{e}{0,T;L^2(\Omega)}\le& \Const{el}\linf{\mathcal{E}(U,g,\mathcal{T})}{0,T}
+2\Const{el}\int_0^T\mathcal{E}(U_t,g_t,\mathcal{T})\\
&+\sqrt{2}\Const{el}\mathcal{E}(U(0),g(0),\mathcal{T})\\
& +\sqrt{2}\ltwo{u_0-U(0)}{}+ C_{a,T} \ltwo{u_1- U_t(0)}{}.
\end{aligned}
\end{equation}
\end{corollary}
\begin{proof} Using (\ref{crucial_ell_rec_property}), $\ltwo{\epsilon}{}$ and $\ltwo{\epsilon_t}{}$
 can be bounded from above using (\ref{apost_bound_bound}).
\end{proof}

\begin{remark}
A bound of the form (\ref{apost_bound_bound}) is only required to
hold for Corollary \ref{apost_semi_thm_wave_l2_cor} to be
valid. Therefore, other available a posteriori bounds for elliptic
problems can be also used; see, e.g.,~\citet{verfurth,ainsworth_oden} and the references therein.
\end{remark}

\section{A posteriori error bounds for the fully discrete problem}\label{full_discrete}

The analysis of \S\ref{abs_apost} is now extended to the case of a
fully-discrete implicit scheme with the aid of a novel three point
space-time reconstruction, satisfying a crucial \emph{vanishing
moment property} in the time variable.

%in Definition~\ref{elliptic_recon_fd}, which is tailored to fit, as
%best as we possibly can, the second derivative structure of the wave
%equation.  An error relation is then derived in
%Proposition~\ref{energy_arg_wave_fd} and the crucial \emph{vanishing
%  moment property} of the space-time reconstruction observed.  We then
%define the indicators and state the main result of this paper:
%Theorem~\ref{apost_fd_thm_wave_l2}.  The proof is then spread across
%several sections.

\begin{definition}[Space-time reconstruction]\label{elliptic_recon_fd}
Let $U^n$, $n=0,\dots,N$, be the fully discrete solution computed by
the method (\ref{fem_fd}), $\Pi^n:L^2(\Omega)\to V_h^n$ be the
orthogonal $L^2$-projection, and $A^n:V_h^n\to V_h^n$ to be the
discrete operator defined by
\begin{equation}
  \label{ell_rec2_fd} \text{for}\ q\in
  V_h^n,\quad \langle A^n q, \chi \rangle = a(q,\chi)\quad \forall
  \chi\in V_h^n.
\end{equation}
We define the \emph{elliptic reconstruction}
$w^n\in H^1_0(\Omega)$, of $U^n$ to be the solution of the elliptic
problem \begin{equation}\label{ell_rec_fd} a(w^n, v) = \langle g^n, v\rangle\quad
\forall v\in H^1_0(\Omega), \end{equation} with
\begin{equation}
  g^n:= A^n U^n-\Pi^n f^n+\bar{f}^n,
\end{equation}
where $\bar{f}^0(\cdot):=f(0,\cdot)$
%$\bar{f}^1(\cdot):=k_1^{-1}\int_{0}^{t^1}\big(f(t,\cdot)-\ha (f(0,\cdot)+A^0 U^0)\big) \ud t$
and
$\bar{f}^n(\cdot):=k_n^{-1}\int_{t^{n-1}}^{t^n}f(t,\cdot) \ud t$ for $n=1,\dots,N$.
Finally, we need to define the \emph{elliptic reconstruction} $\partial w^0\in H^1_0(\Omega)$, of $V^0$ to be
the solution of the elliptic problem
\begin{equation}\label{ell_rec_fd_neq1}
a(\partial w^0, v) = \langle \partial g^0, v\rangle\quad \forall v\in H^1_0(\Omega),
\end{equation}
with
\begin{equation}
\partial g^0:= A^0 V^0-\Pi^0 f^0+f^0.
\end{equation}
\end{definition}

The \emph{time-reconstruction} $\Uh:[0,T]\times \Omega\to \mathbb{R}$ of $\{U^n\}_{n=0}^N$, is defined by
\begin{equation}\label{time_rec}
\Uh(t):= \frac{t-t^{n-1}}{k_n} U^{n} + \frac{t^{n}-t}{k_n} U^{n-1}
-\frac{(t-t^{n-1})(t^{n}-t)^2}{k_n}\partial^2 U^{n},
\end{equation}
for $t\in (t^{n-1},t^n]$, $n=1,\dots, N$, with $\partial^2 U^{n}$ given in (\ref{second_der_def}),
noting that $\partial U^0$ is well defined in (\ref{second_der_def}).
We note that $U$ is a $C^1$-function
in the time variable, with $U(t^n)=U^n$ and
$U_t(t^n)=\partial U^n$ for , $n=0,1,\dots,N$.

We shall also use the time-continuous elliptic reconstruction $\wh$, defined by
\begin{equation}\label{time_elliptic_rec}
\wh(t):= \frac{t-t^{n-1}}{k_n} w^{n} + \frac{t^{n}-t}{k_n} w^{n-1}
-\frac{(t-t^{n-1})(t^{n}-t)^2}{k_n}\partial^2 w^{n},
\end{equation}
noting that $\partial w^0$ is well defined. By construction, this is
also a $C^1$-function in the time variable.

We decompose the error as follows:
\begin{equation}\label{splitting_fd_wave}
e:=\Uh-u=\rho-\epsilon,\ \text{where}\ \epsilon:=\wh-\Uh,\ \text{and}\ \rho:=\wh-u.
\end{equation}

\begin{remark}[Notation overload]
  In this section we use symbols, e.g., $U,w,e,\epsilon,\rho$, that
  where used in \S\ref{abs_apost}, but with a slightly different
  meaning.  Indeed, these are now fully-discrete constructs,
  corresponding in aim and meaning, but different, to their
  semidiscrete counterpart.  It is hoped that this
  overload of notation should not create any confusion.
\end{remark}

\begin{proposition}[Fully-discrete error relation]\label{energy_arg_wave_fd}
For $t\in (t^{n-1},t^{n}]$, $n=1,\dots,N$, we have
\begin{equation}\label{en_ar_wave_fd}
\langle e_{tt},v\rangle +a(\rho,v)=
\langle  (I-\Pi^{n})\Uh_{tt},v\rangle+\mu^n\langle  \partial^2 U^n,\Pi^{n}v\rangle
    +a(\wh- w^{n},v)+\langle   \bar{f}^{n}-f,v\rangle,
\end{equation}
for all $v\in H^1_0(\Omega)$, with $I$ being the identity mapping in $L^2(\Omega)$, and
\begin{equation}
\mu^n(t):=- 6k_n^{-1}(t-t^{n-\ha}),
\end{equation}
 where $t^{n-\ha}:=\ha(t^n+t^{n-1})$.
\end{proposition}
\begin{proof}
Noting that
$\Uh_{tt}(t)= (1+\mu^n(t))\partial^2U^{n}$, for $t\in (t^{n-1},t^{n}]$, $n=1,\dots,N$,
and the identity $a(U^{n},\Pi^{n}v)-\langle \Pi^n f^n-\bar{f}^n, v\rangle = a(w^{n},v)$,
we deduce
\begin{equation}
  \label{rho-term_wave_fd}
  \begin{aligned}
    \langle&  e_{tt},v\rangle +a(\rho,v)
    =\langle  \Uh_{tt},v\rangle +a(\wh,v) - \langle f,v\rangle,\\
    &=\langle  (I-\Pi^{n})\Uh_{tt},v\rangle+\langle  \Uh_{tt},\Pi^{n}v\rangle +a(\wh,v) - \langle f,v\rangle,\\
    &=\langle  (I-\Pi^{n})\Uh_{tt},v\rangle +\mu^n(t)\langle  \partial^2 U^n,\Pi^{n}v\rangle
    \\
    &\phantom=
- a(U^{n},\Pi^{n}v)
+a(\wh,v)+\langle  \Pi^{n}f^{n}-f,v\rangle\\
&=\langle  (I-\Pi^{n})\Uh_{tt},v\rangle+\mu^n(t)\langle  \partial^2 U^n,\Pi^{n}v\rangle
+a(\wh- w^{n},v)+\langle   \bar{f}^{n}-f,v\rangle.
  \end{aligned}
\end{equation}
\end{proof}

\begin{remark}[Vanishing moment property]
The particular form of the remainder $\mu^n(t)$ satisfies the vanishing moment property
\begin{equation}\label{vanishing_moment}
\int_{t^{n-1}}^{t^n} \mu^n(t)\, \ud t =0,
\end{equation}
which appears to be of crucial importance for the optimality of the
a posteriori bounds presented below.
\end{remark}

\begin{definition}[A posteriori error indicators]
  \label{def:aposteriori-error-indicators}
  We define in a list form the error indicators which will form error
  estimator the fully discrete bounds in Theorem \ref{apost_fd_thm_wave_l2}.
  \begin{itemize}
  \item mesh change indicator
    $\eta_1(\tau):=\eta_{1,1}(\tau)+\eta_{1,2}(\tau)$, with
    \begin{equation}
      \eta_{1,1}(\tau):=
      \suint\ltwo{(I-\Pi^{j})\Uh_{t}}{}
      +\mtau \ltwo{(I-\Pi^{m})\Uh_{t}}{},
    \end{equation}
    and
    \begin{equation}
      \eta_{1,2}(\tau):=
      \sum_{j=1}^{m-1}(\tau-t^j)\ltwo{(\Pi^{j+1}-\Pi^{j})\partial U^j}{}
      + \tau\ltwo{(I-\Pi^{0})V^0(0)}{},
    \end{equation}
  \item evolution error indicator
    \begin{equation}
      \displaystyle \eta_2(\tau):=\int_0^{\tau} \ltwo{\mathcal{G}}{},
    \end{equation}
    where $\mathcal{G}:(0,T]\to\mathbb{R}$ with $\mathcal{G}|_{(t^{j-1},t^j]}:=\mathcal{G}^j$,
        $j=1,\dots,N$ and
        \begin{equation}
          \mathcal{G}^j(t):=
          \frac{(t^j-t)^2}{2}\partial g^j-\Big(\frac{(t^j-t)^4}{4k_j}-\frac{(t^j-t)^3}{3}\Big)\partial^2 g^j-\gamma_j,
        \end{equation}
        with $g^j$ as in Definition \ref{elliptic_recon_fd} and
        $\gamma_j:=\gamma_{j-1}+\frac{k_j^2}{2}\partial
        g^j+\frac{k_j^3}{12}\partial^2 g^j$, $j=1,\dots,N$, with
        $\gamma_0=0$;
      \item data error indicator
        \begin{equation}
          \eta_3(\tau)
          :=\frac{1}{2\pi}\sum_{j=1}^{m-1}\Big(\int_{t^{j-1}}^{t^j} k_j^3\ltwo{\bar{f}^{j}-f}{}^2\Big)^{1/2}
          +
          \Big(\mtau k_m^{3}\ltwo{ \bar{f}^{m}-f}{}^2\Big)^{1/2};
        \end{equation}
      \item  time reconstruction error indicator
        \begin{equation}
          \eta_4(\tau)
          :=\frac{1}{2\pi}\sum_{j=1}^{m-1}\Big(\int_{t^{j-1}}^{t^j} k_j^3\ltwo{\mu^j\partial^2 U^j}{}^2\Big)^{1/2}
          +
          \Big(\mtau k_m^{3}\ltwo{\mu^m\partial^2 U^m}{}^2\Big)^{1/2}.
        \end{equation}
  \end{itemize}
\end{definition}

\begin{theorem}[Abstract fully-discrete error bound]\label{apost_fd_thm_wave_l2}
Recalling the notation of Definition~\ref{elliptic_recon_fd} and the
indicators of Definition~\ref{def:aposteriori-error-indicators} we
have the bound
\begin{equation}\label{apost_fd_wave}
\begin{aligned}
\linf{e}{0,t^N;L^2(\Omega)}\le& \linf{\epsilon}{0,t^N;L^2(\Omega)}+\sqrt{2}\Big(\ltwo{u_0-U(0)}{}
+\ltwo{\epsilon(0)}{}\Big)\\
&+2\Big(\int_{0}^{t^N}\ltwo{\epsilon_t}{}+\sum_{i=1}^4\eta_i(t^N)\Big) + C_{a,N} \ltwo{u_1- V^0}{},
\end{aligned}
\end{equation}
where $C_{a,N}:=\min\{2t^N,\sqrt{2C_{\Omega}/\alpha_{\rm min}}\}$.
\end{theorem}
The proof of Theorem \ref{apost_fd_thm_wave_l2} is the content of the remaining of this section.

  Next we set $v=\tilde{v}$ in (\ref{en_ar_wave_fd}) with $\tilde{v}$
  defined by (\ref{tilde_v}) where $\rho$ is defined as in
  (\ref{splitting_fd_wave}) (i.e., the fully discrete $\rho$), assuming that  $t^{m-1}< \tau \le t^m$
  for some integer $m$ with $1\le m \le N$. We integrate the resulting equation with respect to $t$
  between $0$ and $\tau$, to arrive to
  \begin{equation}\label{energy_identity2}
    \int_0^{\tau}\langle e_{tt},\tilde{v}\rangle +\int_0^{\tau}a(\rho,\tilde{v})
    =
    \mathcal{I}_1(\tau) +\mathcal{I}_2(\tau)+\mathcal{I}_3(\tau)+\mathcal{I}_4(\tau),
  \end{equation}
  where
  \begin{equation}
    \label{eqn:proof:fully-discrete:four-components}
    \begin{aligned}
      \mathcal{I}_1(\tau)
      &
      :=\suint\langle  (I-\Pi^{j})\Uh_{tt},\tilde{v}\rangle
      +\mtau\langle  (I-\Pi^{m})\Uh_{tt},\tilde{v}\rangle,
      \\
      \mathcal{I}_2(\tau)
      &
      :=\suint a(\wh- w^{j},\tilde{v})
      +\mtau a(\wh- w^{m},\tilde{v})
      \\
      \mathcal{I}_3(\tau)
      &
      :=\suint \langle  \bar{f}^{j}-f,\tilde{v}\rangle
      +\mtau \langle   \bar{f}^{m}-f,\tilde{v}\rangle,
      \\
      \mathcal{I}_4(\tau)
      &
      :=\suint \mu^j\langle  \partial^2 U^j,\Pi^{j}\tilde{v}\rangle
      +\mtau \mu^m\langle  \partial^2 U^m,\Pi^{m}\tilde{v}\rangle.
    \end{aligned}
  \end{equation}
  In Lemmas \ref{mesh_change_lemma}, \ref{time_estimate_lemma}, \ref{forcing_bound},
  and \ref{time_leftover} we will derive bounds of the form
  \begin{equation}\label{bounds_lemmas}
    \mathcal{I}_i(\tau)
    \le \eta_i(\tau)\max_{0\le t\le T}\ltwo{\rho(t)}{},
  \end{equation}
  for $i=1,2,3,4$.  With the help of these, we will conclude the proof Theorem \ref{apost_fd_thm_wave_l2} at the end of this section.

\begin{lemma}[Mesh change error estimate]\label{mesh_change_lemma}
Under the assumptions of Theorem~\ref{apost_fd_thm_wave_l2} and with
the notation~(\ref{eqn:proof:fully-discrete:four-components}) we have
\begin{equation}
\label{en_ar_wave_fd_first}
\mathcal{I}_1(\tau)
\le \eta_1(\tau)\max_{0\le t\le T}\ltwo{\rho(t)}{}.
\end{equation}
\end{lemma}
\begin{proof}  Observing that the projections $\Pi^j$, $j=1,\dots, N$,
commute with time-differentiation, we integrate by parts with respect to $t$, arriving to
\begin{equation}\label{etaena}
\begin{aligned}
\mathcal{I}_1(\tau)
=&
\suint\langle  (I-\Pi^{j})\Uh_{t},\rho\rangle+
\mtau\langle  (I-\Pi^{m})\Uh_{t},\rho\rangle\\
&+\sum_{j=1}^{m-1}\langle(\Pi^{j+1}-\Pi^{j})\Uh_{t}(t^j),\tilde{v}(t^{j})\rangle
-\langle  (I-\Pi^{0})\Uh_{t}(0),v(0)\rangle.
\end{aligned}
\end{equation}
The first two terms on the right-hand side of (\ref{etaena}) are bounded by
\begin{equation}
\max_{0\le t\le T}\ltwo{\rho(t)}{}
\bigg(
\suint\ltwo{(I-\Pi^{j})\Uh_{t}}{}
+\mtau \ltwo{(I-\Pi^{m})\Uh_{t}}{}
\bigg).
\end{equation}
Recalling the definition of $\tilde{v}$ and that $\Uh(t^j)=\partial U^j$, $j=0,1,\dots, N$,
we can bound the last two terms on the right-hand side of (\ref{etaena}) by
\begin{equation}
\max_{0\le t\le T}\ltwo{\rho(t)}{}
\bigg(
\sum_{j=1}^{m-1}(\tau-t^j)\ltwo{(\Pi^{j+1}-\Pi^{j})\partial U^j}{}
+ \tau\ltwo{(I-\Pi^{0})V^0(0)}{}
\bigg).
\end{equation}
\end{proof}

\begin{lemma}[Evolution error bound]\label{time_estimate_lemma}
Under the assumptions of Theorem~\ref{apost_fd_thm_wave_l2} and with
the notation~(\ref{eqn:proof:fully-discrete:four-components}) we have
\begin{equation}
\label{en_ar_wave_fd_second}
\mathcal{I}_2(\tau)
\le
\eta_2(\tau) \max_{0\le t\le T}\ltwo{\rho(t)}{}.
\end{equation}
\end{lemma}
\begin{proof} First, we observe the identity
\begin{equation}
\wh-w^{j} =
-(t^j-t)\partial w^j+\Big(k_j^{-1}(t^j-t)^3-(t^j-t)^2\Big)\partial^2 w^j,
\end{equation}
on each $(t^{j-1},t^j]$, $j=2,\dots,m$. Hence, from Definition \ref{elliptic_recon_fd}, we deduce
\begin{equation}\label{en_ar_wave_fd_second_one}
a(\wh- w^{j},\tilde{v})=\langle
-(t^j-t)\partial g^j+\Big(k_j^{-1}(t^j-t)^3-(t^j-t)^2\Big)\partial^2 g^j, \tilde{v}\rangle
\end{equation}

The integral of the first component in the inner product on the right-hand side of
(\ref{en_ar_wave_fd_second_one}) with respect to $t$
between $(t^{j-1},t^j]$ is then given by $\mathcal{G}$.
The choice of constants in $\mathcal{G}$ implies that $\mathcal{G}$ is continuous on
$t^j$, $j=1,2,\dots,N$ and $\mathcal{G}(0)=0$.

Hence, integrating by parts on each interval $(t^{j-1},t^j]$, $j=1,\dots,m$,
we obtain
\begin{equation}
\label{en_ar_wave_fd_second_two}
\mathcal{I}_2(\tau)
=\int_0^{\tau} \langle \mathcal{G},\rho\rangle,
\end{equation}
which already implies the result.
\end{proof}

\begin{lemma}[Data approximation error bound]\label{forcing_bound}
Under the assumptions of Theorem~\ref{apost_fd_thm_wave_l2} and with
the notation~(\ref{eqn:proof:fully-discrete:four-components}) we have
\begin{equation}
\label{en_ar_wave_fd_third}
\mathcal{I}_3(\tau)
\le \eta_3(\tau)\max_{0\le t\le T}\ltwo{\rho(t)}{}.
\end{equation}
\end{lemma}
\begin{proof} We begin by observing that
\begin{equation}\label{vanishing_property}
\int_{t^{j-1}}^{t^j}( \bar{f}^{j}-f) =0,
\end{equation}
for all $j=1,\dots,m-1$.
Hence, we have
\begin{equation}
\suint \langle   \bar{f}^{j}-f,\tilde{v}\rangle
=
\suint \langle   \bar{f}^{j}-f,\tilde{v}- \bar{\tilde{v}}^j\rangle,
\end{equation}
where $\bar{\tilde{v}}^j(\cdot):=k_j^{-1}\int_{t^{j-1}}^{t^j}\tilde{v}(t,\cdot) \ud t$.
Using the inequality
%Fubini's Theorem and one-dimensional Poincar\'e inequality (also known as Wirtinger inequality \citet{hardy})
%yield
\begin{equation}
\int_{t^{j-1}}^{t^j}\ltwo{\tilde{v}- \bar{\tilde{v}}^j}{}^2\le \frac{k_j^2}{4\pi^2}
\int_{t^{j-1}}^{t^j}\ltwo{\tilde{v}_t}{}^2,
\end{equation}
and recalling that $\tilde{v}_t=\rho$, we have, respectively,
\begin{equation}
\begin{aligned}
\suint \langle   \bar{f}^{j}-f,\tilde{v}\rangle
&\le
\sum_{j=1}^{m-1}\Big(\int_{t^{j-1}}^{t^j} \ltwo{ \bar{f}^{j}-f}{}^2\Big)^{1/2}
\Big(\int_{t^{j-1}}^{t^j}\ltwo{\tilde{v}- \bar{\tilde{v}}^j}{}^2\Big)^{1/2}\\
&\le\frac{1}{2\pi}
\sum_{j=1}^{m-1}\Big(\int_{t^{j-1}}^{t^j} \ltwo{ \bar{f}^{j}-f}{}^2\Big)^{1/2}
\Big(\int_{t^{j-1}}^{t^j}k_j^2\ltwo{\rho}{}^2\Big)^{1/2}\\
&\le\frac{1}{2\pi}
\sum_{j=1}^{m-1}\Big(\int_{t^{j-1}}^{t^j}
 k_j^3\ltwo{ \bar{f}^{j}-f}{}^2\Big)^{1/2}\max_{0\le t\le T}\ltwo{\rho(t)}{}.
\end{aligned}
\end{equation}
For the remaining term in $\mathcal{I}_3$, we first observe that
\begin{equation}\label{babis}
\mtau \ltwo{\tilde{v}}{}^2\ud t\le \mtau k_m\int_t^{\tau}\ltwo{\rho}{}^2\ud s\ud t\le
k_m^3\max_{0\le s\le T}\ltwo{\rho(t)}{}^2,
\end{equation}
which implies
\begin{equation}
\mtau \langle   \bar{f}^{m}-f,\tilde{v}\rangle
\le
\Big(\mtau k_m^{3}\ltwo{ \bar{f}^{m}-f}{}^2\Big)^{1/2}\max_{0\le t\le T}\ltwo{\rho(t)}{}.
\end{equation}
Recalling $\eta_3$ from
Definition~\ref{def:aposteriori-error-indicators} we conclude the
proof.
\end{proof}

\begin{remark}[The order of the data approximation indicator]\label{f_higher_order}
The choice of the particular combination of functions involving the right-hand side data $f$
in the definition of $g^n$ in the elliptic reconstruction, results to the property (\ref{vanishing_property}).
When $f$ is differentiable, we have $\eta_3(\tau)=O(k^2)$ as $k:=\max_{1\leq j\leq m}k_j\to 0$, and the convergence
is of second order with respect to the maximum time-step. In this case, $\eta_3$ is, therefore, a higher order term.
\end{remark}

\begin{lemma}[Time-reconstruction error bound]\label{time_leftover}
Under the assumptions of Theorem~\ref{apost_fd_thm_wave_l2} and with
the notation~(\ref{eqn:proof:fully-discrete:four-components}) we have
\begin{equation}\label{en_ar_wave_fd_fourth}
\mathcal{I}_4(\tau)
\le \eta_4(\tau)\max_{0\le t\le T}\ltwo{\rho(t)}{}.
\end{equation}
\end{lemma}
\begin{proof}
The method of bounding $\mathcal{I}_4(\tau)$ is similar to that
of Lemma \ref{forcing_bound}, so we shall only highlight the differences.

Recalling the vanishing moment property (\ref{vanishing_moment}) and noting that
$\partial^2U^j$ is piecewise constant in time, we have
\begin{equation}
\suint \mu^j\langle  \partial^2 U^j,\Pi^{j}\tilde{v}\rangle
=
\suint \mu^j\langle  \partial^2 U^j,\Pi^{j}(\tilde{v} -\bar{\tilde{v}}^j)\rangle,
\end{equation}
where $\bar{\tilde{v}}^j(\cdot)=k_j^{-1}\int_{t^{j-1}}^{t^j}\tilde{v}(t,\cdot) \ud t$. Hence,
since $\Pi^j$ commutes with time integration, we obtain
\begin{equation}
\begin{aligned}
\suint \mu^j\langle  \partial^2 U^j,\Pi^{j}(\tilde{v} -\bar{\tilde{v}}^j)\rangle
&\le\frac{1}{2\pi}
\sum_{j=1}^{m-1}\Big(\int_{t^{j-1}}^{t^j} \ltwo{\mu^j\partial^2 U^j }{}^2\Big)^{1/2}
\Big(\int_{t^{j-1}}^{t^j}k_j^2\ltwo{\Pi^j\rho}{}^2\Big)^{1/2}\\
&\le\frac{1}{2\pi}
\sum_{j=1}^{m-1}\Big(\int_{t^{j-1}}^{t^j}
 k_j^3\ltwo{\mu^j\partial^2 U^j }{}^2\Big)^{1/2}\max_{0\le t\le T}\ltwo{\rho(t)}{}.
\end{aligned}
\end{equation}
For the remaining term in $\mathcal{I}_4$, upon using an argument similar to (\ref{babis}), we have
\begin{equation}
\mtau \langle   \mu^m\partial^2 U^m,\Pi^m\tilde{v}\rangle
\le
\Big(\mtau k_m^{3}\ltwo{ \mu^m\partial^2 U^m}{}^2\Big)^{1/2}\max_{0\le t\le T}\ltwo{\rho(t)}{}.
\end{equation}
Recalling the definition of $\eta_4$ in \S\ref{def:aposteriori-error-indicators} we conclude.
\end{proof}

\label{sec:concluding-proof-of-fully-discrete}
 Starting from (\ref{energy_identity2}),
integrating by parts the first term on the left-hand side,
and using the properties of $\tilde{v}$, we arrive at
\begin{equation}
 \int_0^{\tau}\ha \frac{d}{dt}\ltwo{\rho}{}^2
- \int_0^{\tau} \ha \frac{d}{dt}a(\tilde{v},\tilde{v})
= \int_0^{\tau}\langle \epsilon_t,\rho\rangle
+ \langle e_t( 0  ),\tilde{v}( 0  )\rangle+\sum_{i=1}^4\mathcal{I}_i(\tau),
\end{equation}
which implies
\begin{equation}
 \ha\ltwo{\rho(\tau)}{}^2 - \ha\ltwo{\rho(0)}{}^2
+\ha a(\tilde{v}(0),\tilde{v}(0))
= \int_0^{\tau}\langle \epsilon_t,\rho\rangle
+ \langle e_t( 0  ),\tilde{v}( 0  )\rangle+\sum_{i=1}^4\mathcal{I}_i(\tau).
\end{equation}
Hence, we deduce
\begin{multline}
  \label{abstract_bound_fd_mid}
  \ha\ltwo{\rho(\tau)}{}^2 - \ha\ltwo{\rho(0)}{}^2+\ha a(\tilde{v}(0),\tilde{v}(0))
  \\
  \le \max_{0\le t\le T} \ltwo{\rho(t)}{} \Big(\int_0^{\tau}\ltwo{\epsilon_t}{}+\sum_{i=1}^4\eta_i(\tau)\Big)
  + \ltwo{e_t( 0  )}{}\ltwo{\tilde{v}( 0  )}{}.
\end{multline}
We select $\tau=\hat{\tau}$ such that $\ltwo{\rho(\hat{\tau})}{}=\max_{0\le t\le t^N}\ltwo{\rho(t)}{}$.
First, observing that
$\ltwo{\tilde{v}(0)}{}\le \tau\ltwo{\rho(\hat{\tau})}{}$, gives
\begin{equation}
 \frac{1}{4}\ltwo{\rho(\tau)}{}^2 - \ha\ltwo{\rho(0)}{}^2
\le \Big(\int_0^{\tau}\ltwo{\epsilon_t}{}+\sum_{i=1}^4\eta_i(\tau)
+ \tau\ltwo{e_t( 0  )}{}\Big)^2.
\end{equation}
Using the bound $\ltwo{\rho(0)}{}\le \ltwo{e(0)}{}+\ltwo{\epsilon(0)}{}$ and observing that
$e(0)=\hat{U}(0)-u(0)=U^0-u_0$ and that $e_t(0)=\hat{U}_t(0)-u_t(0)=V^0-u_1$, we arrive to
\begin{equation}\label{apost_fd_wave_first}
\begin{aligned}
\linf{e}{0,t^N;L^2(\Omega)}\le& \linf{\epsilon}{0,t^N;L^2(\Omega)}+\sqrt{2}\Big(\ltwo{u_0-U^0}{}
+\ltwo{\epsilon(0)}{}\Big)\\
&+2\Big(\int_0^{t^N}\ltwo{\epsilon_t}{} +\sum_{i=1}^4\eta_i(t^N)+ t^N \ltwo{u_1-V^0}{}\Big).
\end{aligned}
\end{equation}
The second way is completely analogous to the proof of the semidiscrete case. Combining the bounds above suffices to conclude the proof of Theorem \eqref{apost_fd_thm_wave_l2}.

\section{Fully-discrete a posteriori estimates of residual type}
\label{sec:residual-estimates} To arrive to a practical a posteriori
bound for the fully-discrete scheme from
Theorem~\ref{apost_fd_thm_wave_l2}, the quantities involving the
elliptic error $\epsilon$ should be estimated in an a posteriori
fashion: this is the content of Lemmas \ref{epsilon_lemma} and
\ref{epsilon_t_lemma} below, when residual-type a posteriori
estimates are used.

\begin{lemma}[Estimation of the elliptic error]\label{epsilon_lemma}
With the notation introduced in Definition~\ref{elliptic_recon_fd}, we have
\begin{equation}\label{elliptic_error_linf}
\linf{\epsilon}{0,t^N;L^2(\Omega)}+\sqrt{2}\ltwo{\epsilon(0)}{}
\le
\delta_1(t^N)+\sqrt{2}\Const{el}\mathcal{E}^0,
\end{equation}
where
\begin{equation}\begin{aligned}
    \delta_1(t^N)&
    :=\max\Big\{
    \frac{8k_1}{27}\Const{el}\mathcal{E}(V^0,\partial g^0,\mathcal{T}^0),
    \\
    &\Big(\frac{35}{27}+\frac{31}{27}\max_{1\leq j\leq N}\frac{k_j}{k_{j-1}}\Big)
    \max_{0\leq j\leq N}\big(\Const{el}\mathcal{E}^j+
    C_{\Omega}\alpha_{\min}^{-1}\ltwo{\bar{f}^j-f^j}{}\big)
    \Big\},
 \end{aligned}
\end{equation}
with $\mathcal{E}^j:=\mathcal{E}(U^j,A^j U^j-\Pi^j f^j+f^j,\mathcal{T}^j)$, $j=0,1,\dots,N$.
\end{lemma}
\begin{proof} For $t\in(t^{j-1},t^j]$, $j=1,\dots,N$, we have
\begin{equation}
\epsilon(t)=\frac{t-t^{j-1}}{k_j} (w^j-U^{j}) + \frac{t^{j}-t}{k_j} (w^{j-1}-U^{j-1})
-\frac{(t-t^{j-1})(t^{j}-t)^2}{k_j}(\partial^2 w^j-\partial^2U^{j}),
\end{equation}
from which, we can deduce
\begin{equation}
\ltwo{\epsilon(t)}{}\le\max\Big\{\Big(\frac{35}{27}+\frac{31}{27}\max_{1\leq j\leq N}\frac{k_j}{k_{j-1}}\Big)
\max_{0\leq j\leq N}\ltwo{w^j-U^{j}}{}, \frac{8k_1}{27}\ltwo{\partial w^0-V^0}{}\Big\},
\end{equation}
noting that
\begin{equation}
\max_{t\in(t^{j-1},t^j]}\frac{(t-t^{j-1})(t^{j}-t)^2}{k_j}=\frac{4k_j^2}{27}.
\end{equation}
It remains to estimate the terms $\ltwo{w^j-U^{j}}{}$ and $\ltwo{\partial w^0-V^0}{}$. To this end,
recalling the notation of Definition \ref{elliptic_recon_fd}, we define
$w^j_*\in H^1_0(\Omega)$ to be
the solution of the elliptic problem
\begin{equation}
\label{ell_rec_fd_star}
a(w^j_*, v) = \langle A^j U^j-\Pi^j f^j+f^j, v\rangle\quad \forall v\in H^1_0(\Omega),
\end{equation}
for $j=0,1,\dots,N$. Note that, due to the fact that $\bar{f}^0=f^0$,
we have $w^0_*=w^0$.  By construction, we have $a(w^j_*,V)=\langle A^j
U^j-\Pi^j f^j+f^j, V\rangle=a(U^j,V)$ for all $V\in V^j_h$,
$j=0,1,\dots,N$. Hence, $U^j$ is the finite element solution (in
$V^j_h$) of the elliptic boundary-value problem
(\ref{ell_rec_fd_star}).  In view of Theorem
\ref{mraposterioritheorem}, this implies that
\begin{equation}\label{ell_rec_nth_step}
\ltwo{w^j_*-U^j}{}\le \Const{el}\mathcal{E}^j,
\end{equation}
for $j=0,\dots,N$. Similarly, by
construction, we have $a(\partial w^0,V)=\langle A^0 V^0-\Pi^0 f^0+f^0, V\rangle=a(V^0,V)$ for all $V\in V^0_h$.
Hence,
\begin{equation}\label{first_term_ell_rec_bound}
\ltwo{\partial w^0-\partial U^0}{}\le \Const{el}\mathcal{E}(V^0,\partial g^0,\mathcal{T}^0).
\end{equation}
Moreover, since $w^j-w^j_*$ is the solution of an elliptic problem with right hand-side $\bar{f}^j-f^j,$
standard elliptic stability results yield
%
%Let $z,z_* \in H^1_0(\Omega)$ be the solutions to the elliptic problems:
%\begin{equation}\label{ell_prob2}
%-\nabla\cdot(a\nabla z)= r\quad\text{and}\quad -\nabla\cdot(a\nabla z_*)= r_*
%\end{equation}
%with homogeneous Dirichlet boundary conditions, respectively. Then, we have
%\begin{equation}\label{pde_stability}
%\norm{z-z_*}{} \le C_{\Omega}\alpha_{\min}^{-1}\norm{r-r_*}{}.
%\end{equation}
%Theorem \ref{pde_stability_theorem} implies
\begin{equation}\label{pde_stab_estimate}
\ltwo{w^j-w^j_*}{}\le C_{\Omega}\alpha_{\min}^{-1}\ltwo{\bar{f}^j-f^j}{},
\end{equation}
for $j=1,\dots,N$. Finally, using the triangle inequality
\begin{equation}
\ltwo{w^j-U^j}{}\le \ltwo{w^j-w^j_*}{}+\ltwo{w^j_*-U^j}{},
\end{equation}
along with the bounds (\ref{pde_stab_estimate}), (\ref{first_term_ell_rec_bound}) and (\ref{ell_rec_nth_step}),
already implies the result.
\end{proof}

\begin{remark}
The bound (\ref{elliptic_error_linf}) contains both the \emph{elliptic estimators} $\mathcal{E}(\cdot,\cdot,\cdot)$
and the data-oscillation terms $\ltwo{\bar{f}^j-f^j}{}$ which are, in general,
of first order with respect to the time-step.
The data-oscillation terms are expected to dominate the data error indicator $\eta_3$ (cf. Remark \ref{f_higher_order}).
On the other hand, if the numerical scheme (\ref{fem_fd})
is altered so that $f^j=\bar{f}^j$ (as done, e.g., in \citet{baker_wave}),
then the data-oscillation terms in (\ref{elliptic_error_linf}) vanish.
Similar remarks apply to the result of Lemma 4.12 below.
\end{remark}

For each $n=1,\dots, N$, we denote by $\hat{\mathcal{T}}^{n}$ the finest
common coarsening of $\mesh{n}$ and $\mesh{n-1}$, and by $\hat{V}_h^n:=V_h^n\cap V_h^{n-1}$, the
corresponding finite element space, along with the orthogonal $L^2$-projection
operator $\hat{\Pi}^n:L^2(\Omega)\to\hat{V}_h^n$.

\begin{lemma}[Estimation of the time derivative of the elliptic error]\label{epsilon_t_lemma}
With the notation introduced in \S\ref{elliptic_recon_fd} we have
\begin{equation}\label{elliptic_error_linftwo}
\int_0^{t^N}\ltwo{\epsilon_t}{}
\le
\delta_2(t^N),
\end{equation}
where
\begin{equation}\begin{aligned}
\delta_2(t^N):=\frac{2}{3}\sum_{j=0}^{N}(2k_j+k_{j+1})
\Big( \Const{el}\mathcal{E}^j_{\partial}+C_{\Omega}\alpha_{\min}^{-1}\ltwo{\partial f^j-\partial \bar{f}^j}{}\Big),
 \end{aligned}
\end{equation}
with
\begin{equation}\mathcal{E}_{\partial}^j
:=\mathcal{E}(\partial U^j,\partial (A^j U^j)-\partial(\Pi^j f^j)
+ \partial f^j,\hat{\mathcal{T}}^j),\quad  j=0,1,\dots,N.
\end{equation}
\end{lemma}
\begin{proof} For $t\in(t^{j-1},t^j]$, $j=1,\dots,N$, we have
\begin{equation}
\epsilon_t=\partial w^j-\partial U^{j}
-k_j^{-1}(t^{j}-t)(t^j-2t^{j-1}+t)(\partial^2 w^j-\partial^2U^{j}),
\end{equation}
from which, we deduce
\begin{equation}\label{epsilon_t_bounds}
\int_{t^{j-1}}^{t^j}
\ltwo{\epsilon_t}{}\le \frac{4k_j}{3}\ltwo{\partial w^j-\partial U^{j}}{}
+\frac{2k_j}{3}\ltwo{\partial w^{j-1}-\partial U^{j-1}}{},
\end{equation}
noting that
\begin{equation}
\int_{t^{j-1}}^{t^j}k_j^{-2}(t^{j}-t)(t^j-2t^{j-1}+t) \ud t=\frac{2k_j}{3}.
\end{equation}
Combining (\ref{epsilon_t_bounds}) for $j=1,\dots,N$, we arrive to
\begin{equation}\label{epsilon_t_bounds_all}
\int_0^{t^N}
\ltwo{\epsilon_t}{}\le
\frac{2}{3}\sum_{j=0}^{N}(2k_j+k_{j+1})\ltwo{\partial w^j-\partial U^{j}}{},
\end{equation}
with $k_0=0$ and $k_{N+1}=0$.

It remains to estimate the terms $\ltwo{\partial w^j-\partial U^{j}}{}$. To this end,
recalling the definition of the functions $w^j_*\in H^1_0(\Omega)$ from
the proof of Lemma \ref{epsilon_lemma} and, since $\hat{V}_h^j:=V_h^j\cap V_h^{j-1}$, we have
$a(w^j_*,V)=a(U^j,V)$ for all $V\in \hat{V}^j_h$ and $a(w^{j-1}_*,V)=a(U^{j-1},V)$ for all $V\in \hat{V}^j_h$
, for $j=1,\dots,N$. Therefore, we deduce
\begin{equation}
a(\partial w^j_*,V)=a(\partial U^j,V)\quad\text{for all}\quad V\in \hat{V}^j_h,
\end{equation}
for $j=1,\dots,N$, i.e., $\partial U^j$ is the finite element solution in $\hat{V}^j_h$ of the boundary-value
problem
\begin{equation}
a(\partial w^j_*,V)=\langle \partial (A^j U^j)-\partial(\Pi^j f^j)
+ \partial f^j, v\rangle\quad \forall v\in H^1_0(\Omega).
\end{equation}
In view of Theorem \ref{mraposterioritheorem}, this implies that
\begin{equation}\label{ell_rec_nth_step_two}
\ltwo{\partial w^j_*-\partial U^j}{}\le \Const{el}\mathcal{E}_{\partial}^j,
\end{equation}
for $j=1,\dots,N$. We also recall that, by construction, we have
$a(\partial w^0,V)=a(V^0,V)$ for all $V\in V^0_h$. Hence, (\ref{first_term_ell_rec_bound}) also holds.

Moreover, since
\begin{equation}
a(\partial w^j,V)=\langle \partial (A^j U^j)-\partial(\Pi^j f^j)
+ \partial \bar{f}^j, v\rangle\quad \forall v\in H^1_0(\Omega),
\end{equation}
$j=1,\dots,N$, (cf. Definition \ref{elliptic_recon_fd}).
%Theorem \ref{pde_stability_theorem} implies
As in (\ref{pde_stab_estimate}), elliptic stability implies
\begin{equation}\label{pde_stab_estimate_two}
\ltwo{\partial w^j-\partial w^j_*}{}\le
C_{\Omega}\alpha_{\min}^{-1}\ltwo{\partial \bar{f}^j-\partial f^j}{},
\end{equation}
for $j=1,\dots,N$ and, using the triangle inequality
\begin{equation}
\ltwo{\partial w^j-\partial U^j}{}\le \ltwo{\partial w^j-\partial w^j_*}{}+\ltwo{\partial w^j_*-\partial U^j}{},
\end{equation}
along with the bounds (\ref{pde_stab_estimate_two}), (\ref{first_term_ell_rec_bound}) and (\ref{ell_rec_nth_step_two}),
already implies the result.
\end{proof}

\begin{theorem}[Fully-discrete residual-type a posteriori bound]\label{apost_fd_thm_wave_l2_practical}
With the same hypotheses and notation as
in Theorems \ref{apost_fd_thm_wave_l2} and \ref{mraposterioritheorem}, we have the bound
\begin{equation}\label{apost_fd_wave_practical}
\begin{aligned}
\linf{e}{0,t^N;L^2(\Omega)}\le& \delta_1(t^N)+\sqrt{2}\Const{el}\mathcal{E}^0+\sqrt{2}\ltwo{u_0-U(0)}{}
\\
&+2\delta_2(t^N)+2\sum_{i=1}^4\eta_i(t^N) + C_{a,N} \ltwo{u_1- V^0}{},
\end{aligned}
\end{equation}
where $\delta_1, \mathcal{E}^0$  are defined in Lemma \ref{epsilon_lemma}, $\delta_2$ is
defined in Lemma \ref{epsilon_t_lemma}, and $\eta_i$, $i=1,2,3,4$ after (41) respectively.
\end{theorem}
\begin{proof}
  Combining Theorem \ref{apost_fd_thm_wave_l2} with the bounds derived for $\epsilon$
  in Lemma~\ref{epsilon_lemma}, and $\epsilon_t$ in Lemma~\ref{epsilon_t_lemma},
  we arrive to an a posteriori error bound.
\end{proof}

\section{Final remarks}\label{conclusions}
The design and implementation of adaptive algorithms for the wave
equation based on rigorous a posteriori error estimators is a
largely unexplored subject, despite the importance of these problems
in the modelling of a number of physical phenomena. To this end, this
work presents rigorous a posteriori error bounds in the
$L^{\infty}(L^2)$-norm for second order linear hyperbolic
initial/boundary value problems. The use of a novel space-time-reconstruction technique, that hinges on the one-field formulation of the problem, appears to be generic and it is expected to be applicable to second order hyperbolic problems with, posssibly nonlinear spatial operators or other with spatial discretisations. Although the case of residual-type estimators has been demonstrated above, it is evident that Theorem \eqref{apost_fd_thm_wave_l2} can be combined with a variety of 
other a posteriori estimators for elliptic problems. The derived bounds appear to be of optimal order, although no efficiency bounds are presented; this would be an interesting direction of further research. It is worth noting, however, some of the terms appearing in the a posteriori bound presented above are completely analogous to some of the terms in the respective a posteriori bounds from \citet{bernardi_suli}, which are, in turn, shown to be efficient.  The numerical implementation of the proposed bounds
in the context of adaptive algorithm design for second order
hyperbolic problems deserves special
attention and will be considered elsewhere.

\section*{Acknowledgements}
EHG. acknowledges the support of the Nuffield Foundation, UK, and
of the Foundation for Research and Technology-Hellas, Heraklion,
Greece. OL acknowledges the partial support of the Royal Society UK and of the
Foundation for Research and Technology-Hellas, Heraklion, Greece,
where the initial steps of this work were made. CM acknowledges the support of the London Mathematical Society, Universities of Leicester and Sussex, UK, and supported in part by
the European Union grant No. MEST-CT-2005-021122.

%%%%%%%%%%%%%%%%%%%%%%%%%%%%%%%%%%%%%%%%%%%%%%%%%%%%%%%%%%%%%%%%%%%%%%%%
\bibliographystyle{IMANUM-BIB}

\end{document}